\documentclass{article}
\usepackage[preprint, nonatbib]{neurips_2020}

\usepackage{amsmath,amsthm,amsfonts,amssymb,amscd}
\usepackage{mathrsfs,bm}

\usepackage{hyperref}
\usepackage{lastpage}
\usepackage{enumerate}
\usepackage{fancyhdr}
\usepackage{mathrsfs}
\usepackage{xcolor}
\usepackage{graphicx}
\usepackage{listings}
\usepackage{setspace}
\usepackage{algorithm, algorithmic}
\usepackage{braket}

\newcommand{\slfrac}[2]{\left.#1\middle/#2\right.}

\theoremstyle{definition}
\newtheorem{definition}{Definition} 
\newtheorem{assumption}{Assumption} 

\DeclareMathAlphabet\mathbfcal{OMS}{cmsy}{b}{n}

\usepackage{array, cite}
\usepackage{float}
\usepackage{multirow}
\usepackage{multicol}
\usepackage[caption=false,font=footnotesize]{subfig}

\newcommand{\bDelta}{\bm{\Delta}}
\newcommand{\bLambda}{\bm{\Lambda}}

\newcommand{\bPsi}{\bm{\Psi}}

\newcommand{\ointl}{\oint\limits}

\newcommand{\bdelta}{\bm{\delta}}

\newcommand{\tlbG}{\tilde{\bG}}
\newcommand{\tlbX}{\tilde{\bX}}

\newcommand{\tlbg}{\tilde{\bg}}

\newcommand{\bA}{\mathbf{A}}
\newcommand{\bD}{\mathbf{D}}

\newcommand{\bF}{\mathbf{F}}
\newcommand{\bG}{\mathbf{G}}
\newcommand{\bH}{\mathbf{H}}
\newcommand{\bI}{\mathbf{I}}

\newcommand{\bL}{\mathbf{L}}

\newcommand{\bO}{\mathbf{O}}
\newcommand{\bP}{\mathbf{P}}
\newcommand{\bQ}{\mathbf{Q}}
\newcommand{\bR}{\mathbf{R}}
\newcommand{\bS}{\mathbf{S}}

\newcommand{\bW}{\mathbf{W}}
\newcommand{\bX}{\mathbf{X}}

\newcommand{\bZ}{\mathbf{Z}}

\newcommand{\bU}{\mathbf{U}}
\newcommand{\bC}{\mathbf{C}}

\newcommand{\ba}{\mathbf{a}}

\newcommand{\bc}{\mathbf{c}}
\newcommand{\bd}{\mathbf{d}}

\newcommand{\mbf}{\mathbf{f}}
\newcommand{\bg}{\mathbf{g}}
\newcommand{\bh}{\mathbf{h}}

\newcommand{\bn}{\mathbf{n}}

\newcommand{\bw}{\mathbf{w}}
\newcommand{\bx}{\mathbf{x}}
\newcommand{\by}{\mathbf{y}}
\newcommand{\bz}{\mathbf{z}}

\newcommand{\bbR}{\mathbb{R}}

\newcommand{\calC}{\mathcal{C}}

\newcommand{\calE}{\mathcal{E}}

\newcommand{\calG}{\mathcal{G}}

\newcommand{\calN}{\mathcal{N}}

\newcommand{\calV}{\mathcal{V}}

\newcommand{\calO}{\mathcal{O}}

\newcommand{\barbx}{\bar{\bx}}

\newcommand{\tlbx}{\tilde{\bx}}

\newcommand{\td}{{\text{d}}}

\newcommand{\bzero}{\mathbf{0}}
\newcommand{\bone}{\mathbf{1}}

\newcommand{\suml}{\sum\limits}
\newcommand{\minl}{\min\limits}

\newcommand{\liml}{\lim\limits}

\newcommand{\bigcapl}{\bigcap\limits}

\newcommand{\tr}{\text{Tr}}

\newcommand{\tth}{\text{th}}

\newcommand{\nwl}{\nonumber\\}

\newtheorem{theorem}{Theorem}
\newtheorem{lemma}{Lemma}

\newcommand{\eg}{{{\em e.g.}}}

\newcommand{\nf}[1]{\|#1\|_\mathrm{F}} 
\newcommand{\nt}[1]{\|#1\|_2}

\newcommand{\sbjt}{\mbox{{s.t.}}}

\DeclareMathOperator*{\argmin}{arg\,min}

\newcommand{\brond}{{\mbox{\boldmath $\partial$}}}
\title{Decentralized Constrained Optimization: \\ Double Averaging and Gradient Projection}

\author{
        Firooz~Shahriari-Mehr,~David Bosch~and~Ashkan Panahi\\ 
        Department of Computer Science and Engineering \\
        Chalmers University of Technology \\
        Gothenburg, Sweden \\
        \texttt{Firooz, Davidbos, Ashkan.Panahi@Chalmers.se}
       }

\begin{document}
\maketitle

\begin{abstract}

In this paper, we consider the convex, finite-sum minimization problem with explicit convex constraints over strongly connected directed graphs. The constraint is an intersection of several convex sets each being known to only one node. 
To solve this problem, we propose a novel decentralized projected gradient scheme based on local averaging  and prove its convergence using only local functions' smoothness.
Experimental studies demonstrate the effectiveness of the proposed method in both constrained and unconstrained problems.

\end{abstract}

\section{Introduction}
\label{sec: Intro}


In the past decade, decentralized optimization techniques 
have attracted significant interest ~\cite{Nedic2020,xin2020general}.
In this setting, multiple computing nodes are involved,  and there is no coordinator (central) node with which all nodes communicate. 
A fairly general framework for decentralized optimization problems is given by
\begin{equation}
\label{eq: dec_opt}
    \minl_{\bx \in \bbR^m} \quad f(\bx) \triangleq \suml_{v=1}^{M} f_v(\bx),
\end{equation}
where $M$ is the total number of the nodes in the network, $\bx \in \bbR^m$ is called the global optimization variable, and $f(.)$ is the global objective function which has a finite-sum structure. 
Here, each node $v$ has access to its local function $f_v:\bbR^m \rightarrow \bbR$ and  communicates to its neighbors $\calN_v$ to achieve an \emph{optimal consensus solution}. A natural extension to this setup is when $\bx$ in problem \eqref{eq: dec_opt} is required to lie in an intersection of several convex sets, i.e. $\bx\in \bigcap_{v=1}^N S_v$, and each constraint $S_v$ is known only to one node. Applications of this setup, that we refer to as the Decentralized Constrained Optimization Problem (DCOP), are ubiquitous,  
\eg~smart grid control~\cite{alizadeh2012demand, guan2010energy}, optimal energy management~\cite{nguyen2019distributed}, sensor networks~\cite{bazerque2009distributed}, and support vector machines~\cite{bishop2006pattern}.
However, practical approaches to solving it have not been extensively discussed in the literature. This paper responds to this shortcoming by providing a numerical method to solve DCOP with guaranteed convergence properties.

Node communication is a crucial factor in the design of decentralized optimization techniques, which is represented by either a directed or undirected communication graph. 
Earlier studies on decentralized techniques considered static undirected  graphs, which indicate that each communication link between two nodes is time invariant and bi-directional, meaning that both nodes can send and receive information. This assumption is not compatible with many practical applications, such as broadcast channels with no return link, or communication failures leading to uni-directional links~\cite{Xin2018}. These problems intrigue researchers to propose decentralized methods considering directed graphs as the underlying communication network, where each communication link in the network is uni-directional.  
For simplicity, decentralized methods in which a directed graph represents their node communication are called \emph{decentralized directed methods}, throughout this paper. In the same way, we refer to \emph{decentralized undirected methods}. We address the more general case of decentralized directed scenarios, while the undirected cases follow as a special case.

For undirected communication graphs, efficient optimization techniques with provable convergence properties exist based on suitable iterative averaging over neighbor nodes \cite{boyd2006randomized, nedic2009distributed, shi2015extra}. The averaging procedure is mathematically represented by the so-called gossip matrices, which are compatible with the network structure, double stochastic, and symmetric \cite{shi2015extra, nedic2017achieving}. The required gossip matrices can be constructed using  Laplacian or Metropolis matrices for undirected graphs.
Such gossip matrices are not compatible with directed graphs, as they require asymmetry and finding doubly stochastic matrices is not straightforward, often requiring distributed and iterative numerical procedures such as iterative weight balancing \cite{gharesifard2012distributed}.
For this reason, practical schemes utilize row stochastic or column stochastic matrices, instead of using doubly stochastic matrices. In this case, 
convergence bounds comparable  to the undirected scenarios, even in the absence of constraints, are lacking to the best of our knowledge. We further address this issue by proposing a novel double-averaging scheme, similar to the so-called push-pull approach \cite{pu2020push}, which takes both row and column stochastic matrices into account, and at the same time enjoys superior convergence guarantees.

\subsection{Contributions}

The main contributions of the paper are summarized as follows: 
\begin{itemize}
    \item We propose a novel algorithm, called DAGP,   to solve the problem of decentralized constraint optimization. Our scheme  employs double averaging and projection onto convex sets.
    It extends the tracking approach, which has been proposed for the first time in \cite{shi2015extra,nedic2017achieving}, to constrained problems and can benefit a fixed step size and fast convergence. 
    \item In contrast to the previously proposed methods in the literature, our method  simultaneously considers a directed communication graph and individual constraints at each node. 
    \item We show that our technique is applicable to generic constrained convex problems, lacking strong convexity, while maintaining the convergence rate of order $\calO(\slfrac{1}{\sqrt{n}})$, under mild conditions. We are not aware of any decentralized unconstrained method over directed communication graphs with similar established convergence properties. 
    \item We present experiments for constrained decentralized optimization problems on directed graphs, where DAGP outperforms the existing algorithms. We also conduct experiments on  unconstrained problems, where DAGP performs similarly to the state of the art, decentralized optimization algorithms.
\end{itemize}

\subsection{Literature Review}
In  this section, we review the decentralized methods in the existing literature. We organize our review into three parts: techniques on directed and undirected graphs, and  decentralized constrained methods. 
Several classes of methods exist in the literature that are not in the scope of this paper,~\eg, methods considering time-varying graphs~\cite{nedic2014distributed, nedic2017achieving}, local functions with a finite-sum structure~\cite{Mokhtari2016,Hendrikx2019,Hendrikx2020a,Xin2020}, or compressed communication~\cite{Koloskova2019, beznosikov2020biased}.

\subsubsection{Decentralized optimization over undirected graphs}
The algorithms for undirected communication graphs can be divided into several categories. 
First, the decentralized gradient decent methods including~\cite{nedic2009distributed, lobel2010distributed}  use diminishing step size for convergence to the exact solution of the problem. The diminishing step size leads to practical difficulties with step tuning but establish the convergence rate of $\calO(\slfrac{\log n}{\sqrt{n}})$ in a convex and smooth setting and $\calO(\slfrac{\log n}{n})$ in a strongly convex and smooth setting. 
The second category refers to the methods that use gradient tracking technique and leverage the gradient information at all nodes to estimate the gradient of the global function~\cite{shi2015extra, nedic2017achieving, Qu2018}. These methods use fixed step sizes and achieve linear convergence rate, i.e. $\calO(\mu^{n})$ $\mu < 1,$ in a strongly convex and smooth setting. \cite{Qu2018} has also shown a sublinear rate of convergence, i.e. $\calO(\slfrac{1}{n})$, when the functions are not strongly convex. The dual-based methods~\cite{Seaman2017, uribe2020dual, hendrikx2020optimal} are included in the third group. Although these methods are optimal and have linear convergence rate, they need to compute some computationally costly oracles, \eg~the gradient of a conjugate function, which is not practical in some applications. 

\subsubsection{Decentralized optimization over  directed  graphs}

Earlier methods for directed problems apply the so-called push-sum protocol \cite{kempe2003gossip} to decentralized gradient descent methods to tackle the problem of computing a doubly stochastic gossip matrix for directed graphs~\cite{tsianos2012push, nedic2014distributed}. These methods utilize a column stochastic matrix, but a diminishing step size is still vital for their convergence. The methods based on the push-sum protocol converge with order of $\calO(\slfrac{\log n}{n})$ for smooth and strongly convex functions. To achieve fixed step size, \cite{Xi2018, nedic2017achieving} have put push-sum protocol and gradient tracking technique together and have respectively proposed the DEXTRA and Push-DIGing algorithms. These algorithms achieved a linear rate of convergence in a smooth and strongly convex setting. DEXTRA suffers  theoretical limitations  on the step size, namely a feasible step size might not exist in some cases. \cite{Xi2018} has proposed the ADD-OPT algorithm to solve this problem. This algorithm also enjoys linear convergence for strongly convex functions. 
Recently, methods based on two gossip matrices, one column stochastic and the other one row stochastic, have been proposed in the literature~\cite{pu2020push,xin2018linear}, called Push-Pull methods. These methods also have linear convergence in a smooth and strongly convex setting. 
Our algorithm is similar to push-pull methods as it can be applied with similar underlying matrices.

\subsubsection{Decentralized Constrained optimization} 

Despite extensive studies on decentralized optimization, there exist few papers that consider constraints explicitly. It is worth mentioning that 
a straightforward approach to solving constrained problems is to add the indicator functions of the constraint sets to the problem, then apply the methods proposed for the unconstrained problem. This approach requires methods that are applicable to non-smooth and non-strongly convex functions with unbounded (sub)gradients due to indicator function characteristics. For this reason, we note that utilizing the previously mentioned methods does not guarantee convergence.
\cite{ram2010distributed} is among the first papers incorporating the projection and averaging approaches, but it assumes that the constraint set is identical at all nodes. This leads to a problem when the projection onto the constraint set is not computationally efficient. In response, the projected subgradient algorithm has been proposed, which assumes that the constraint set is different and distributed among all nodes~\cite{nedic2010constrained}. 
This paper is similar in setup to ours, but does not provide a precise convergence rate.  Moreover, convergence of local variables to a consensus stopping point  is proven only in two special cases: when the constraints are identical, or when the graph is fully connected. 
As the constraint at each node might be an intersection of several constraints, or in some applications, the nodes do not have access to all of their local constraints at each iteration, \cite{lee2013distributed} has proposed a randomized projection scheme. This algorithm suffers from the same limitations as~\cite{nedic2010constrained}, i.e. the proof is only reliable for fully connected networks or a setting with identical constraints at each node.
All the above-mentioned methods use a diminishing step size as they do not leverage any gradient tracking technique. Moreover, they assume that the underlying communication graph is undirected. \cite{xi2016distributed} has proposed the DDPS algorithm, which is applicable when the communication graph is directed. However, this algorithm uses diminishing step size as well, and its convergence rate is of order $\calO(\slfrac{\log n}{\sqrt{n }})$.
Moreover, it is subject to the restrictive assumption that the constraints are identical.
%
%
There are also methods with different problem description, such as composite constrained optimization \cite{lu2020computation}. These methods consider undirected communication graphs, and they are different from our problem, in nature.

\subsection{Paper Outline}
The rest of the paper is organized as follows. In the following, some preliminary definitions and notations are introduced. 
The DAGP algorithm is proposed in section~\ref{sec: algorithm}, along with theoretical convergence analysis. The proofs of all Lemmas and Theorems are provided in the Appendix.
Finally, section~\ref{sec: results} is devoted to the numerical studies.






\theoremstyle{definition}
\begin{definition}[Normal cone and Projection Operator]
For a closed convex set $S\subset\bbR^n$, the normal cone of $S$ is given by
\begin{equation*}
\partial I_S(\bx) = 
\begin{cases}
\;\emptyset  & \bx \notin S \\
\left\{\bg \in \bbR^n \vert\;\forall \bz \in S,\;\bg^T(\bz - \bx) \leq 0 \right\} & \bx \in S 
\end{cases}.
\end{equation*}
Moreover, the projection of a vector $\bx \in \bbR^n$ onto $S$ is computed by
\[ 
P_S(\bx) = \argmin\limits_{\by \in S} \nt{\by -\bx}^2.
\]
\end{definition}


\theoremstyle{definition}
\begin{definition}[Graph Theory]
A directed graph is shown by $\calG = (\calV, \calE)$, where $\calV =\{1,\dots,N\}$ is a set of all nodes, and $\calE \subseteq \calV \times \calV$ is a set of ordered pairs of distinct nodes, called edges. 
A directed path between two distinct nodes $u,v \in \calV$ is a sequence of nodes $(u=v_0, v_1, \dots, v_k= v)$ such that each pair $(v_i,v_{i+1})$ is an edge in $\calE$. A
graph $\calG$ is strongly connected if for any two distinct nodes $u, v \in \calV$, there exist a directed path between $u$ and $v$. 
The adjacency matrix denoted by $\bA = [a_{ij}]$ is an asymmetric matrix, where $a_{ij}$ is $+1$ if $(i,j) \in \calE$, and $0$ otherwise. 

In this paper, each pair shows a communication link between two distinct nodes, and $(i,j)$ is a pair in $\calE$, if there is a link from node $j$ to node $i$. With this intuition, the $i$th row of an adjacency matrix shows from which nodes it can receive information, and constitute the incoming neighbors of node $i$ called $\calN_i^{\text{in}} = \left\{j \vert (i,j) \in \calE \right\}$. On the other hand, the $i$th column of an adjacency matrix shows to which nodes it can send information, and constitute the outgoing neighbors of node~$i$ called $\calN_i^{\text{out}} = \left\{ j\vert (j,i)  \in \calE\right\} $. The in-degree and out-degree of node $i$ are defined as the cardinality of $\calN_i^{\text{in}}$ and $\calN_i^{\text{out}}$, respectively. Consequently, two Laplacian matrices can be defined as 
\begin{equation*}
    \bL^{\text{in}} = \bD^{\text{in}} - \bA, \qquad \bL^{\text{out}} = \bD^{\text{out}} - \bA,
\end{equation*}
where $\bD^{\text{in}}$ is the in-degree diagonal matrix; that is $d_{ii}^{\text{in}} = \left\vert\calN_i^{\text{in}}\right\vert$, and $\bD^{\text{out}}$ is
defined in a similar way. 
$\bL^{\text{in}}$ and $\bL^{\text{out}}$ have zero row-sum and column-sum characteristics and their scaled versions are used in this paper. 

\end{definition}


\subsection{Mathematical Notation}
In this paper, bold lowercase  and  uppercase letters are used to respectively represent vectors and matrices. $w_{v u}$ denotes the element at the $v^\tth$ row and the $u^\tth$ column of the matrix $\bW$. $\bW^{T}$ shows the  transpose of $\bW$, and $\ker(\bW)$ is its right null space, meaning that  $\bx \in \ker(\bW)$, if and only if $\bW\bx = \bzero$. 
$\bone_n$ and  $\bzero_n$ respectively denote the $n-$dimensional vectors of all ones and zeros. The index $n$ may be dropped if there is no risk of confusion. Furthermore, $\bO$ denotes a matrix with all zero elements.
The Euclidean inner product of vectors is denoted by $\langle.,.\rangle$. Matrix inner product is  denoted by $\langle\bA,\bC\rangle = \tr(\bA\bC^T) $. 
In this paper, subscript generally defines the iteration number, and superscript defines the node number, \eg, $\nabla f_v(\bx^v_{n})$ indicates the gradient of the node $v$'s local function at its local variable at iteration $n$. Finally,
$\delta_{n,0}$ is the Kronecker delta function.


\section{Problem setting and Proposed Algorithm}
\label{sec: algorithm}
In this section, we propose a new algorithm to solve DCOP, considering directed graphs as a communication network between the nodes.
The proposed algorithm is called \emph{DAGP} due to Double Averaging and Gradient Projection approaches
used in the iterative equations, which are introduced in the subsection \ref{subsec: dagp}.

\subsection{Problem Formulation}
Decentralized Constraint Optimization Problem (DCOP) is formulated as
%
\begin{equation}
\label{eq: dec_constrained}
     \minl_{\bx \in \bbR^{m}}  \;\;  f(\bx) \triangleq \suml_{v=1}^{M} f_v(\bx) \quad
     \sbjt  \quad \bx\in \bigcapl_{v=1}^M S_v,
\end{equation}
where $S_v$ is a closed convex set, and the intersection of all these constraint sets is called the feasible  set.
Note that without loss of generality, 
the number of constraints is equal to the number of functions.
This allows us to assume $M$ nodes in our setting, each having access to one function and one constraint. 
%
Note that a different number of constraints can still be considered in our setup, as each node may further access a composite constraint (i.e. an intersection of simpler constraints), or multiple nodes may share identical constraints (i.e. $S_v=S_u$) or have trivial constraints $S_v=\bbR^m$. Nevertheless, we merely require the projection operator to the constraint set $S_v$ of each node to be available, and neglect further possible structures in them.

In decentralized optimization, each node stores and updates a local variable $\bx^v$ as its solution. The nodes should achieve consensus, i.e. $\bx^v$s must converge to an equal stopping point $\bx^*$, which is further required to be a feasible and optimal solution of \eqref{eq: dec_constrained}.  


\subsection{DAGP Algorithm}
\label{subsec: dagp}
The DAGP algorithm does the following updates in each iteration, $\forall v \in \calV$.
%
\begin{alignat}{3}
    &\bz^v   && = \; &&\bx^v_n - \suml_{u \in \calN_v^{\text{in}}} w_{v u}\bx^u_n - \mu \left( \nabla f_v(\bx^v_n) -  \bg^v_n   \right) \label{eq: z_update}\\
    &\bx_{n+1}^v && = \; &&P_{S_v} \left( \bz^v \right) \label{eq: x_update} \\
    &\bg_{n+1}^v && = \; &&\bg_n^v + \rho \left[ \nabla f_v(\bx^v_n)  - \bg^v_n + \frac{1}{\mu}\left(\bz^v - \bx^v_{n+1} \right)   \right] \nonumber \\ 
    &   &&  && \; + \alpha \left( \bh_n^v - \bg_n^v \right) \label{eq: g_update} \\
    &\bh_{n+1}^v && = \; &&\bh_n^v -\suml_{u \in \calN_v^{\text{in}}} q_{v u}(\bh_n^u - \bg_n^u) \label{eq: h_update}
\end{alignat}
%
To interpret the algorithm, consider
the above update equations. Take into account that the node $v$ has access to $\left(\bx^u_n, \bh^u_n - \bg^u_n\right)$, $\forall u \in  \calN_v^{\text{in}}$, as each node $u$ broadcasts this pair of messages to its out-neighbors.
First weighted averaging happens in \eqref{eq: z_update}, where each node computes a weighted average of its local variable and its in-neighbors' local variables. This averaging is a basis for achieving consensus, and $\bW = \left[ w_{v u} \right]$ must have zero row-sum structure to achieve this objective.
Then, the resulting averaged vector is aligned with the negative of the augmented local descent direction $\left( \nabla f_v(\bx^v_n) -  \bg^v_n   \right)$ scaled by a fixed step size $\mu$.
The resulting solution $\bz^v$ is projected onto the local constraint set in \eqref{eq: x_update}. Therefore from the second iteration, the local variables at each iteration lie in their own local constraint set, but not necessarily in the feasible set of the problem in \eqref{eq: dec_constrained}.
One of the novelties of this paper is the definition of additional variables $\bh^v_n$ together with $\bg^v_n$ to push the algorithm towards an optimal consensus solution. 

Vectors $\bg^v$ act  as the  memory of the algorithm, which preserve and track the previous information of local functions' gradients and \emph{feasible directions}, i.e. $\nabla f_v$ and $\bz^v - \bx^v$, respectively.
To reach an optimal solution, the gradients and feasible directions of all nodes must  be aggregated. 
This is achieved by adding the term $\alpha (\bh^v - \bg^v)$ to  \eqref{eq: g_update}, where $\bh^v$ is updated using \eqref{eq: h_update}. 
$\bh^v$  propagates the information of the gradients and feasible directions of other nodes through the second weighted averaging using $\bQ = \left[ q_{v u} \right]$.
%
In \eqref{eq: h_update}, we further require $\bQ$ to be a matrix with zero column-sum structure. Then, $\sum_{v \in \calV} \bh^v$ will not change over time. 
We also require $\bh^v_0$s be initialized in a way that satisfy $\sum_{v\in \calV} \bh^v = \bzero$. The easiest way is to initialize them with zero vectors. In this way, when $\bg^v$s converges to $\bh^v$s, their summation over all nodes will be equal to zero. This in turn, leads to satisfying the optimality condition of the problem as $\bg^v$s contain gradients and feasible directions of the problem. This is further elaborated in Theorem \ref{thm: consensus_opt}, and Appendix A.

\subsection{Convergence Analysis}
\label{sec: consensus_optimali}
In this section, we discuss the convergence properties of DAGP. We present two results. First in Theorem~\ref{thm: consensus_opt}, we prove that if the iterates of DAGP  converge, any stopping point is an optimal and consensus solution of the problem in \eqref{eq: dec_constrained}. Then, in Theorem~\ref{thm: rate}, we prove the convergence rate of our proposed algorithm in a smooth and convex setting. 

\subsubsection{Assumptions}
 \label{sec: assumptions}
We proceed by formalizing the adopted assumptions 
as follows. 


\begin{assumption}
The nodes will communicate across a strongly connected directed graph $\calG = (\calV,\calE)$.
This assumption guarantees the sufficient information flow between the nodes as there exists a directed path between every two nodes in the graph. As a result, the nodes can achieve consensus. 
\end{assumption}

\begin{assumption}
The optimization problem \eqref{eq: dec_constrained} is feasible and attains a finite optimal value $f^*=f(\bx^*)$ at an optimal feasible solution $\bx^*$ satsfying the optimality condition:
 \[\bzero \in \sum_{v=1}^M \left( \partial I_{S_v}(\bx^*) + \nabla f_v(\bx^*)\right).\] 
\end{assumption}

\begin{assumption}
There are two weight matrices $\bW$ and $\bQ$ with a similar sparsity pattern to the adjacency matrix $\bA$ of $\calG$. They further satisfy the zero row-sum and zero column-sum structure, respectively. 
The first one is required for achieving consensus, and the second one is required for attaining the optimality of the solution.
Moreover, we assume that $\ker(\bQ) = \ker(\bW^T)$ and $\ker(\bW)=\mathrm{span}\{\bone\}$.
\end{assumption}

\begin{assumption}
The functions $f_v(.)$ are convex, differentiable and $L$-smooth.
\end{assumption}

Now, we define the matrices $\bR$ and $ \bP$ as
\begin{equation} \label{eq: RP_define}
\bR = \left[
\begin{array}{cccc}
        \bO & \bO & \bO & \bO   \\
        \bI & \bO & \bO & \bO  \\
        -\frac{\rho}{\mu}\bI & \frac{\rho}{\mu}(\bI-\bW) & \bI & \alpha \bI \\
        \frac{\rho}{\mu}\bI & -\frac{\rho}{\mu}(\bI-\bW) & \bO & (1-\alpha)\bI-\bQ 
\end{array} \right],
\qquad
\bP = \left[
\begin{array}{c}
     \bI  \\
     \bO  \\
     \bO  \\
     \bO 
\end{array}
\right].
\end{equation}
Moreover, for an arbitrary positive value of $\eta$, matrix $\bS$ is computed as
%
%
\begin{equation}\label{eq: S_paper_define}
    \bS=\left[\begin{array}{cccc}
        \left(1-\frac{L\mu}2\right)\bI-M\eta\left(\bI-\frac 1M\bone\bone^T\right) & -\frac 12(\bI-\bW)+\frac {L\mu}2\bI  & -\frac {\mu}2\bI  & \bO\\
        -\frac 12(\bI-\bW^T)+\frac {L\mu}2\bI & -\frac{L\mu}2\bI & \bO &\bO \\
        -\frac{\mu}2\bI & \bO & \bO & \bO \\
         \bO & \bO & \bO& \bO
    \end{array}\right].
\end{equation}
Please check Appendix B for understanding the evolution of these matrices.

\begin{assumption}
There exists a strictly positive constant $C$ such that for every value of $\beta>0$ and in a small neighborhood of $z=0$ on the complex plane, $1$ is not an eigenvalue of the matrix
\begin{equation}
    \left[\begin{array}{ccc}
        \bI\ & \bO\ & \bO   
    \end{array}\right]\bF^{-1}(z,\beta)\left[\begin{array}{c}
          -(C+\beta)\bI \\
            \bI\\
            \bO
    \end{array}\right],
\end{equation}
where $\bF$ is defined as 
\begin{equation} \label{eq: F_define}
    \bF(z,\beta)=\left[\begin{array}{ccc}
        \bS & z^{-1}\bI-\bR^T & \bO\\
        z\bI-\bR & \bO & -\bP \\
        \bO & -\bP^T & -\beta\bI
    \end{array}\right].
\end{equation}

\end{assumption}

\subsubsection{Main Results}
\label{sec:results}
Here, we present the main theoretical results and postpone the details and proofs to the Appendix.
\begin{theorem} 
\label{thm: consensus_opt}
Let Assumptions 3 and 4 hold. 
If the iterates of DAGP algorithm converge, any stopping point is an optimal and consensus solution of the decentralized constrained optimization problem in \eqref{eq: dec_constrained}, i.e. $\bx^v = \bx^*, \; \forall v \in \calV$, and $\bx^*$ satisfies the sufficient optimality conditions. 

\end{theorem}

We also present guarantees for the rate of convergence. 
\begin{theorem}  \label{thm: rate}
Let all the assumptions hold. Define $\barbx^v_N=\frac 1N\suml_{n=0}^{N-1}\bx^v_n$ and $\barbx_N=\frac 1M\suml_v\barbx^v_N$. Then,  $\|\barbx_N-\barbx^v_N\|^2=O(\frac 1N)$, $\mathrm{dist}^2(\barbx_N, S_v)=O(\frac 1N)$ and 
\begin{equation}
 \left|\suml_v f_v(\barbx^v_N)-\suml_v f_v(\bx^*)\right|=O\left(\frac 1{\sqrt{N}}\right).
\end{equation} 
\end{theorem}

\section{Experimental Results}
\label{sec: results}
We evaluate and compare the performance of the DAGP algorithm in two scenarios; decentralized constraint and unconstrained problems.
%
In the first experiment, which contains examples with synthetic data, we consider constraints and examine the convergence and feasibility gap of the DAGP compared to the DDPS algorithm. In the second experiment, we solve the classical logistic regression problem, which is unconstrained. In the latter, we compare our algorithm to the state of the art decentralized unconstrained optimization algorithms over directed graphs, namely the ADD-OPT and Push-Pull algorithms. 
%

\begin{figure}
	\centering
	\subfloat[First setup]{\includegraphics[width = 0.4\linewidth]{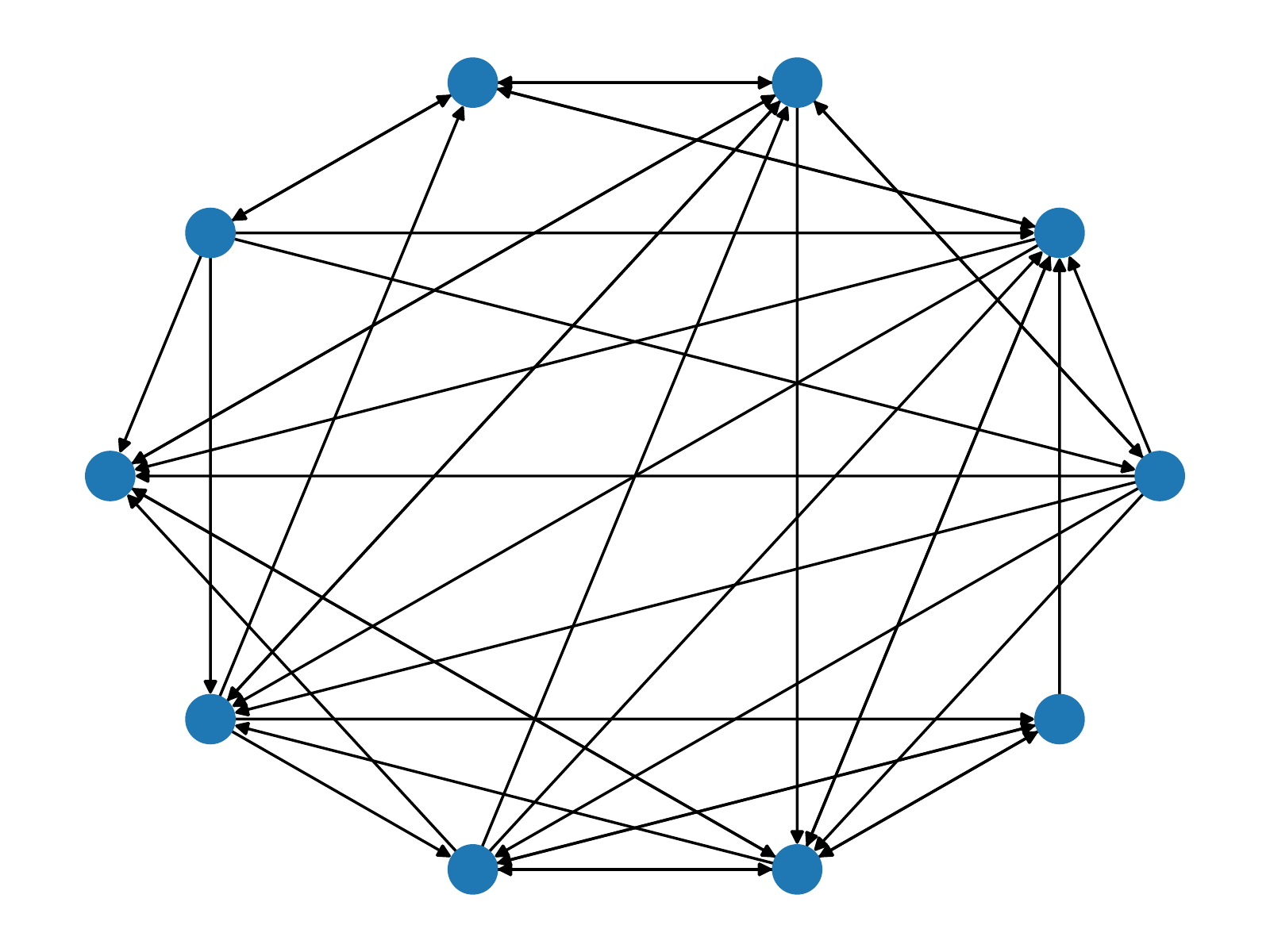} \label{graph1}} \hfill
	\subfloat[Second setup]{\includegraphics[width = 0.4\linewidth]{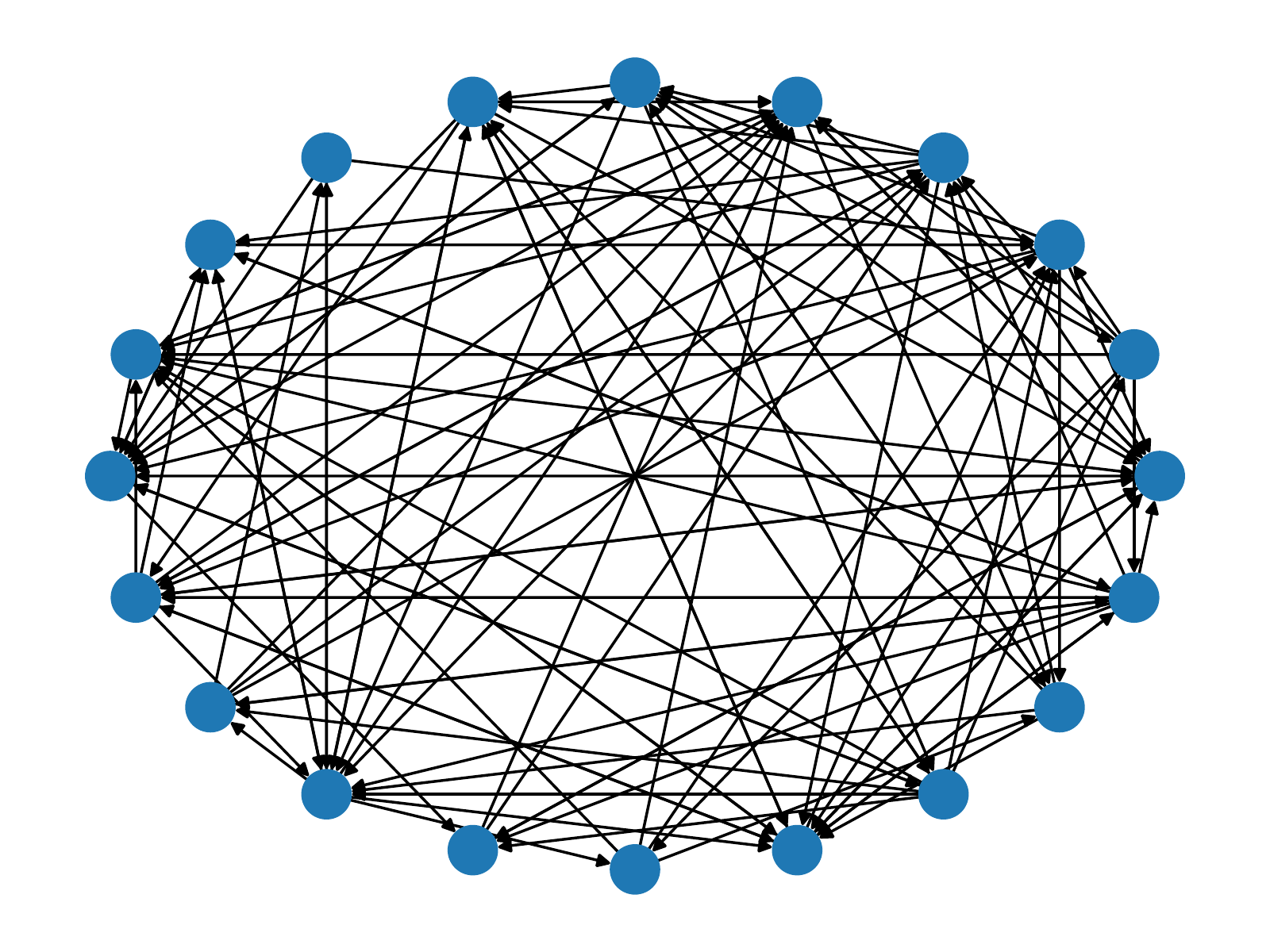} \label{graph2}} \hfill
	
	\subfloat[Logistic Regression]{\includegraphics[width = 0.4\linewidth]{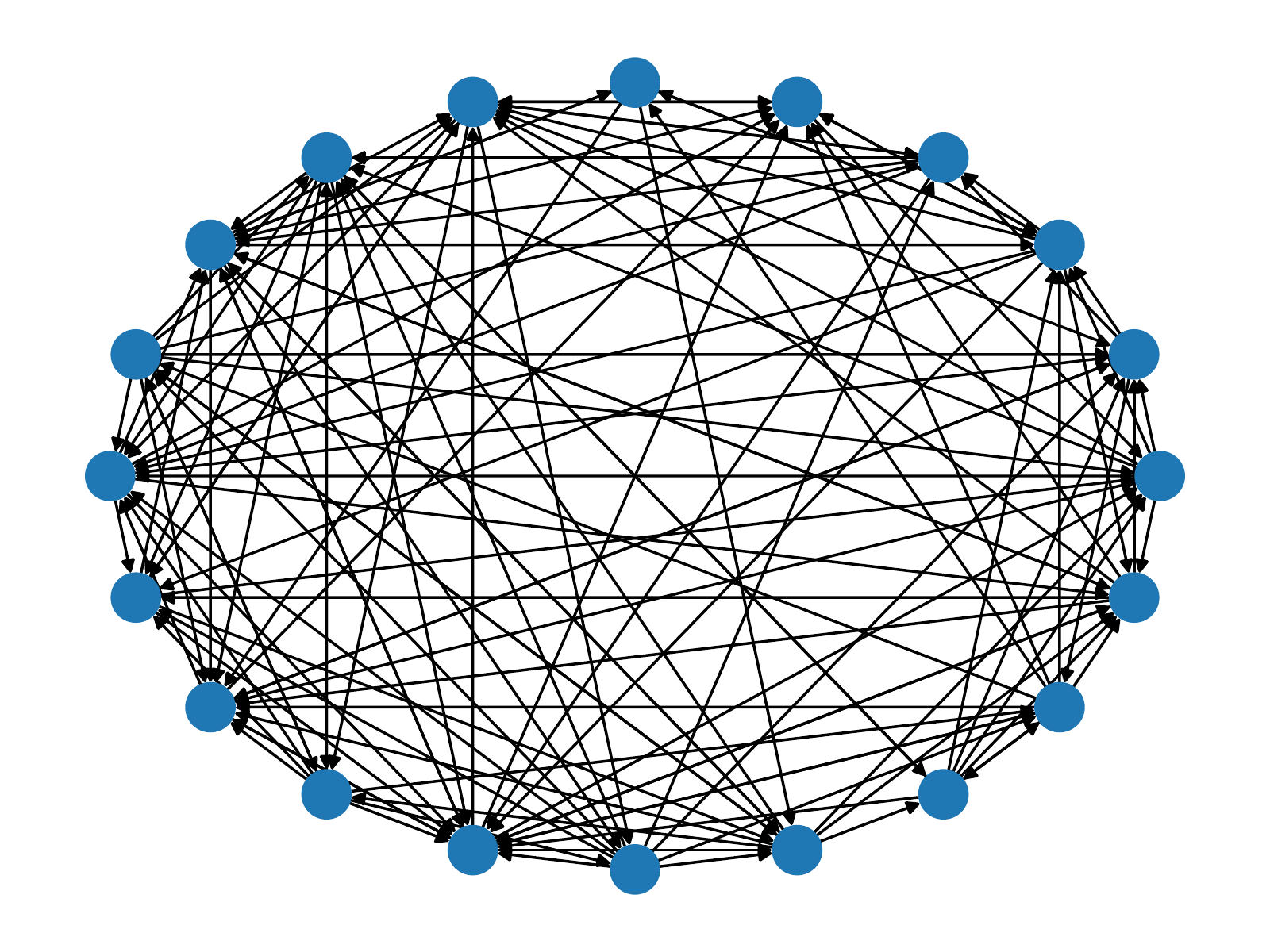}  \label{graph3}}

	\caption{Directed random graphs used in our experiments.}
	\label{fig: graphs}
\end{figure}

In all algorithms used in the first experiment, the parameters  are hand-tuned in a way that the algorithms achieve their best performance, leading to a fair comparison.  In the real-world logistic regression problem, hand-tuning parameters is not computationally feasible due to the size of the experiment. Therefore, an appropriate step size is selected for all algorithms in the second experiment. 
Different algorithms use different matrices for averaging. In this paper, we respectively use 
$\slfrac{\bL^{\text{in}} }{ 2\bd_{\text{max}}^{\text{in}} }$ and $\slfrac{\bL^{\text{out}} }{ 2\bd_{\text{max}}^{\text{out}} }$
as zero row-sum and zero column-sum matrices, where $\bd_{\text{max}}^{\text{in}}$ and $\bd_{\text{max}}^{\text{out}}$ are the largest diagonal elements of $\bL^{\text{{in}}}$ and $\bL^{\text{{out}}}$. By subtracting these matrices from the identity matrix, row stochastic and column stochastic matrices used in this paper  are computed.
Moreover, random directed and strongly connected graphs are used in our experiments, which are shown in Fig~\ref{fig: graphs}. The numerical experiments are described next. We repeated each experiment multiple times, but only one instance from each experiment is presented as the difference between individual runs was minimal.

\subsection{Numerical Results}
In this experiment, which contains two setups,  we consider synthetic functions and constraints. In our setups, there are $M$ nodes, each having access to one function and one constraint. The nodes communicate over a randomly generated graph and their local variables $\bx^v$s are of size $m$, which their initial vectors are generated randomly from zero mean and unit variance normal distribution. The functions are selected to be smooth, but not strongly convex as follows:
\begin{equation}
f_v(\bx) = \log\big(\cosh( \ba_v^T\bx - b_v)\big),
\label{eq: synthetic}
\end{equation}
where $\ba_v$s and $b_v$s are randomly generated from a zero mean and unit variance normal distribution.  
Moreover, we choose randomly generated linear constraints $\bc_v^T\bx - d_v \leq 0$ since their orthogonal projection operator is simple to compute.%
\footnote{
In all simulations, $\bc_v$ and $d_v$ are selected such that their intersection not being an empty set, \eg, $\bc_v$s are generated randomly, then $d_v$s are selected such that $\bc_v^T\bx \leq d_v$ for one arbitrary vector $\bx$. 
}

In the first setup, $m = 20$ and $M = 10$, while in the second one, $m=10$ and $M=20$. These parameters are chosen since the feasible set in the second experiment is significantly smaller than the feasible set in the first experiment. 
In both setups, the objective value and the distance to the feasible set, called feasibility gap, are reported, both computed at $\bar\bx = \sum_{v\in\calV}\bx^v$. 

\begin{figure}[!t]
	\centering
	\subfloat[Objective Value]{\includegraphics[width = 0.49\linewidth]{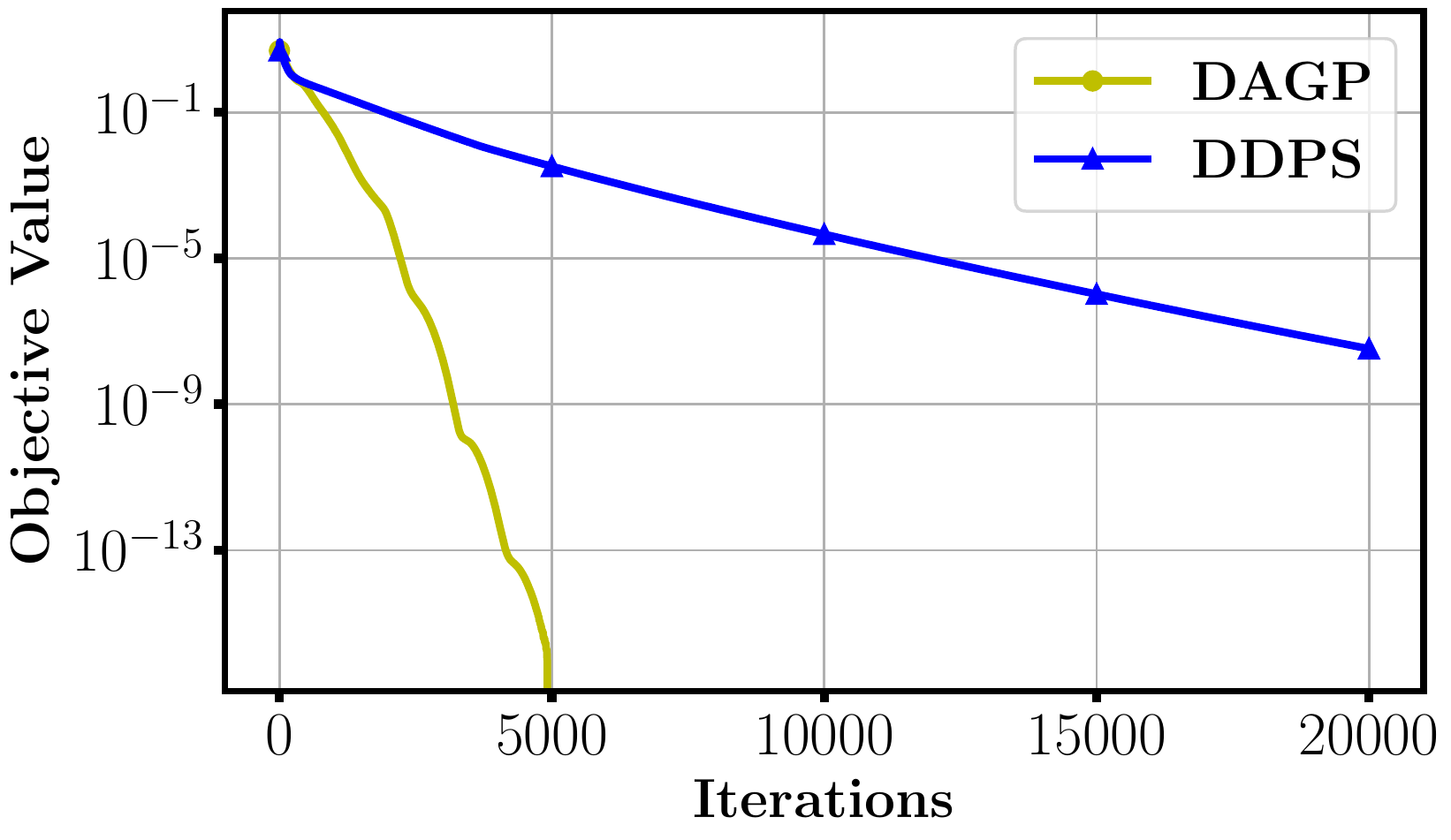} \label{fig: objectiveN10}} 
	\subfloat[Feasibility Gap]{\includegraphics[width = 0.49\linewidth]{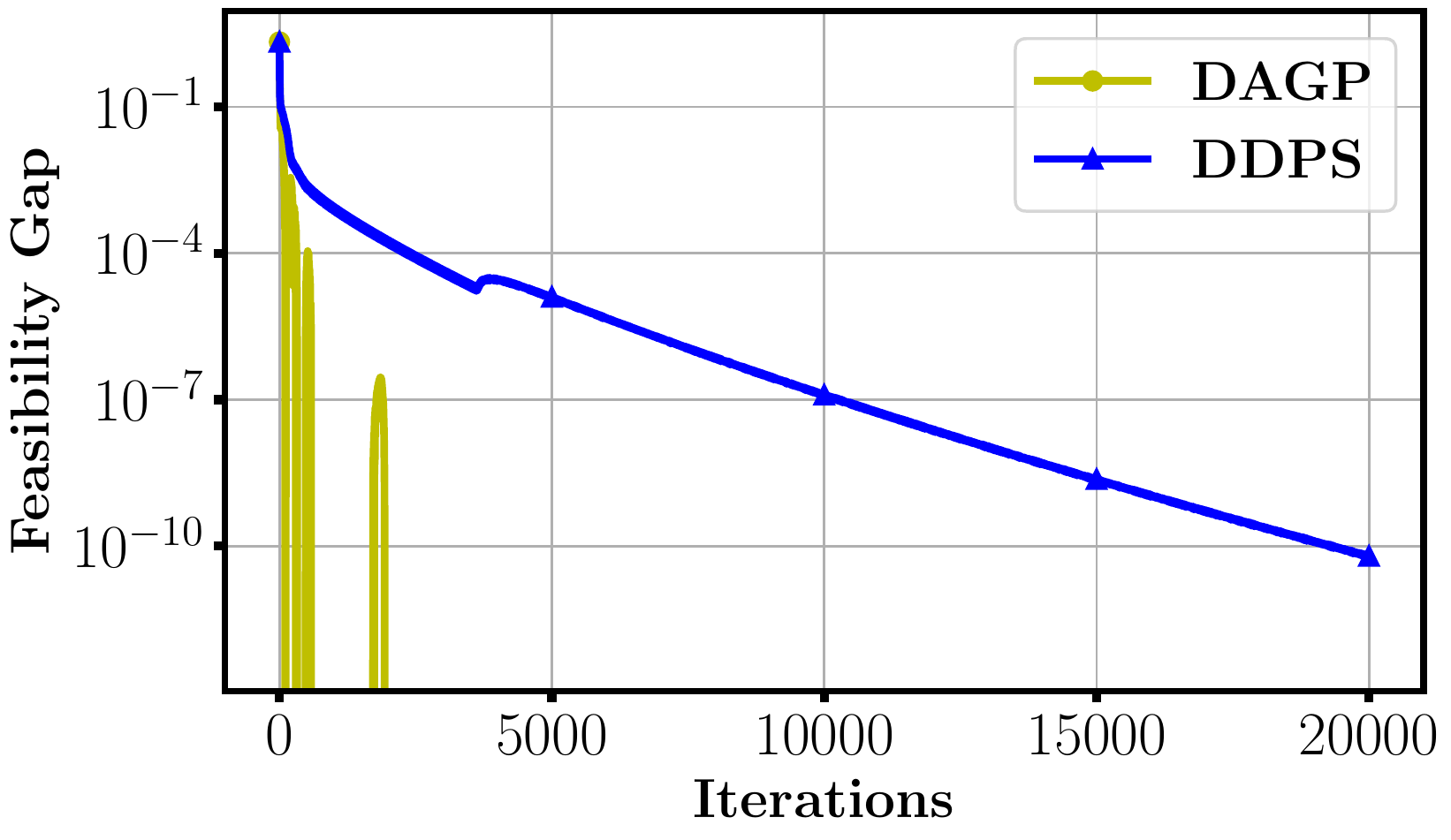}   \label{fig: fesgpN10} }
    
	\subfloat[Consensus solution]{\includegraphics[width = 0.49\linewidth]{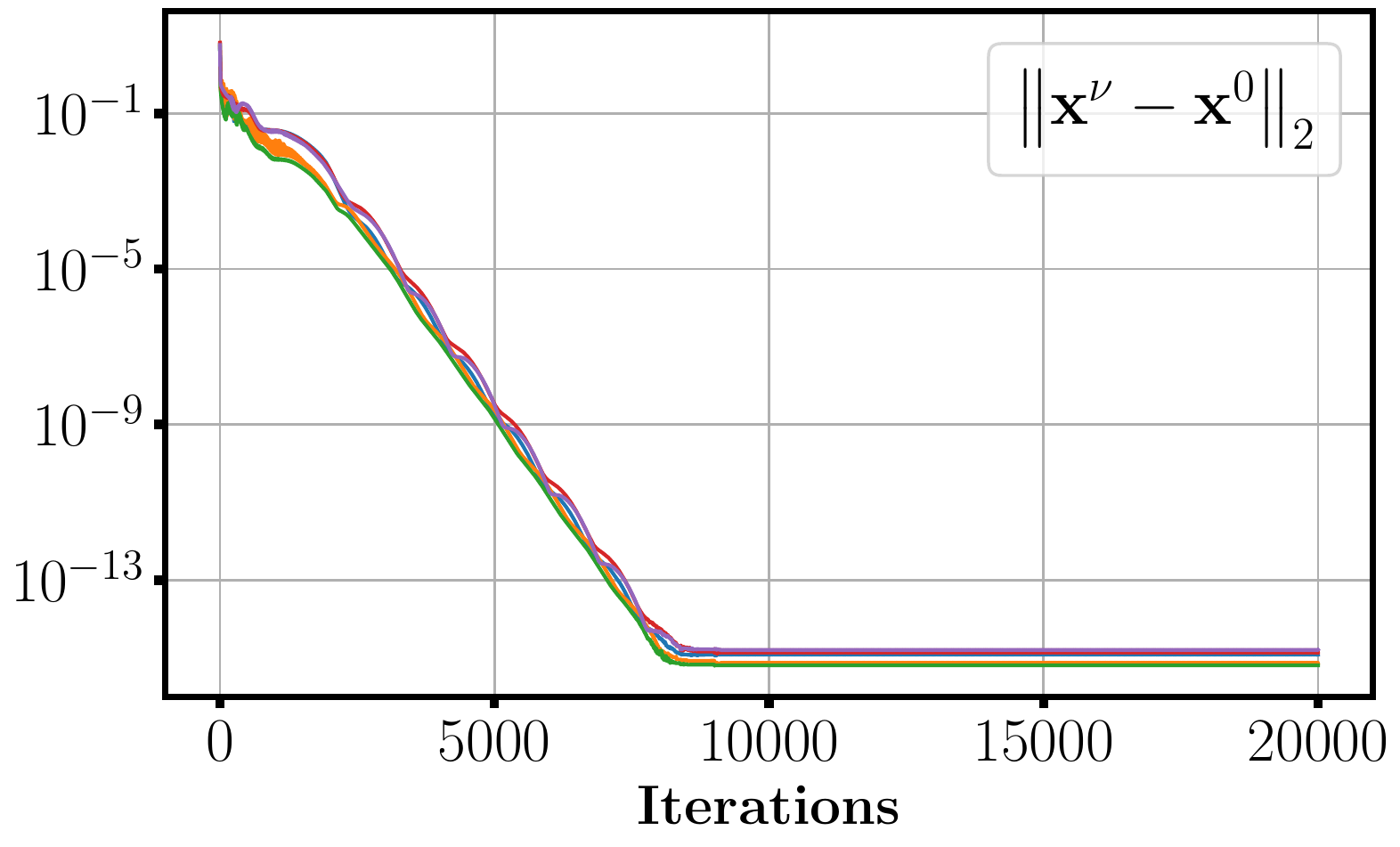}  \label{fig: cons}} 
	\subfloat[Optimal solution]{\includegraphics[width = 0.49\linewidth]{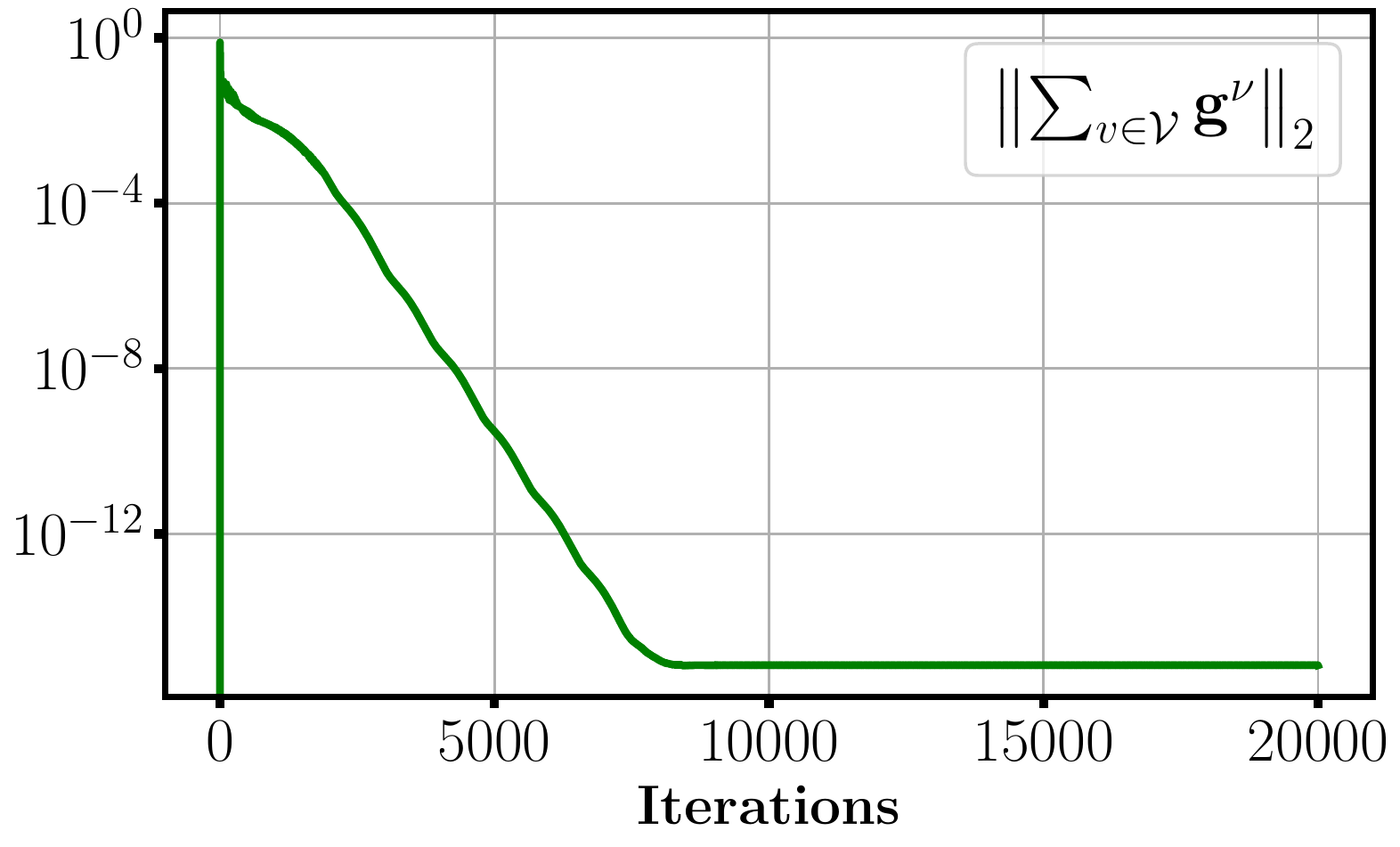} \label{fig: gsum}}
	\caption{First setup results with $m=20$ and $M=10$. Local variables move to a consensus and optimal stopping point in DAGP, while they move to  a sub-optimal point in DDPS.}
	\label{fig:toy1}
\end{figure}

\begin{figure}
    \centering
	\subfloat[Objective Value]{\includegraphics[width = 0.475\linewidth]{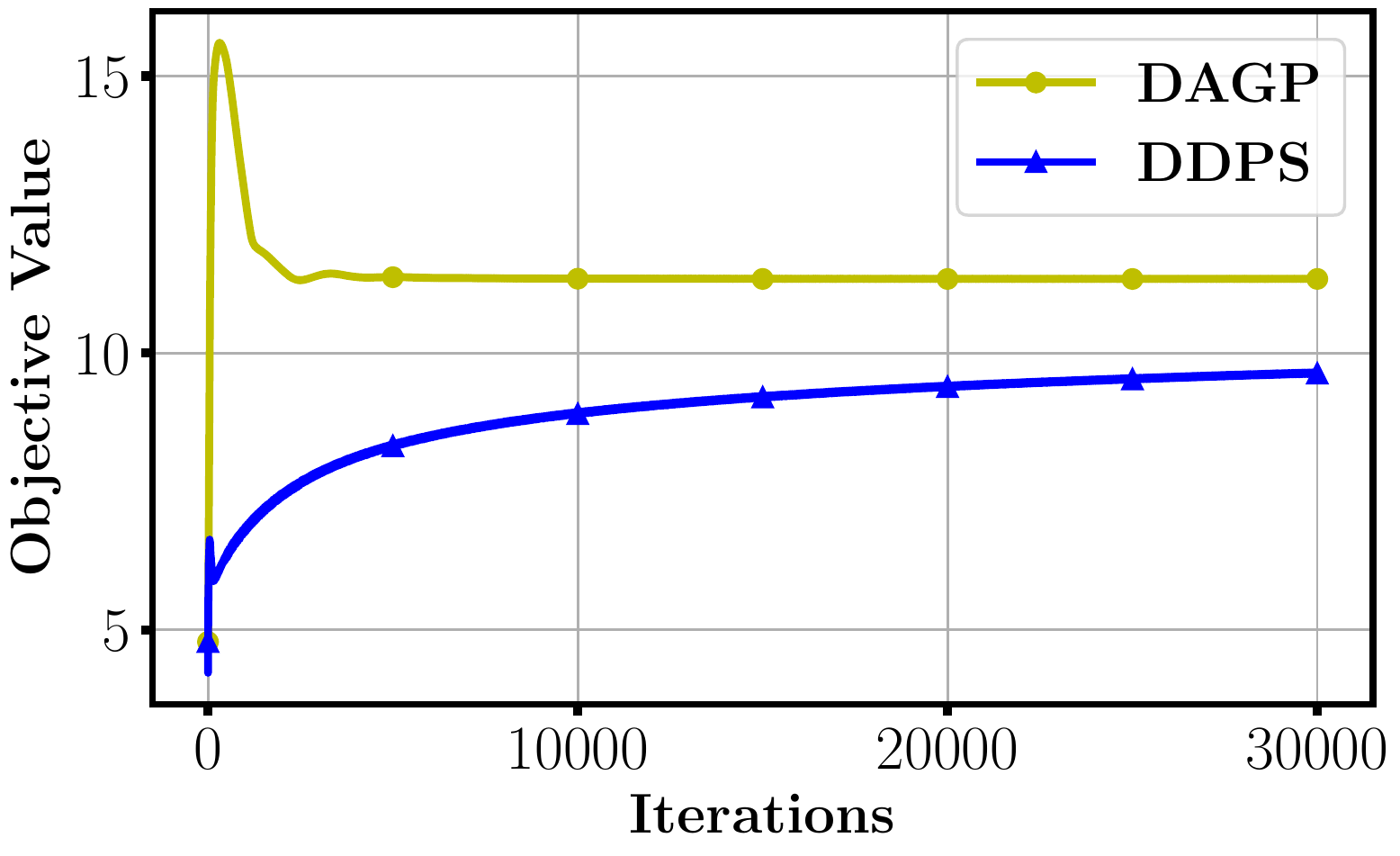} \label{fig: ObjectiveN20}} 
	\subfloat[Feasibility Gap]{\includegraphics[width = 0.49\linewidth]{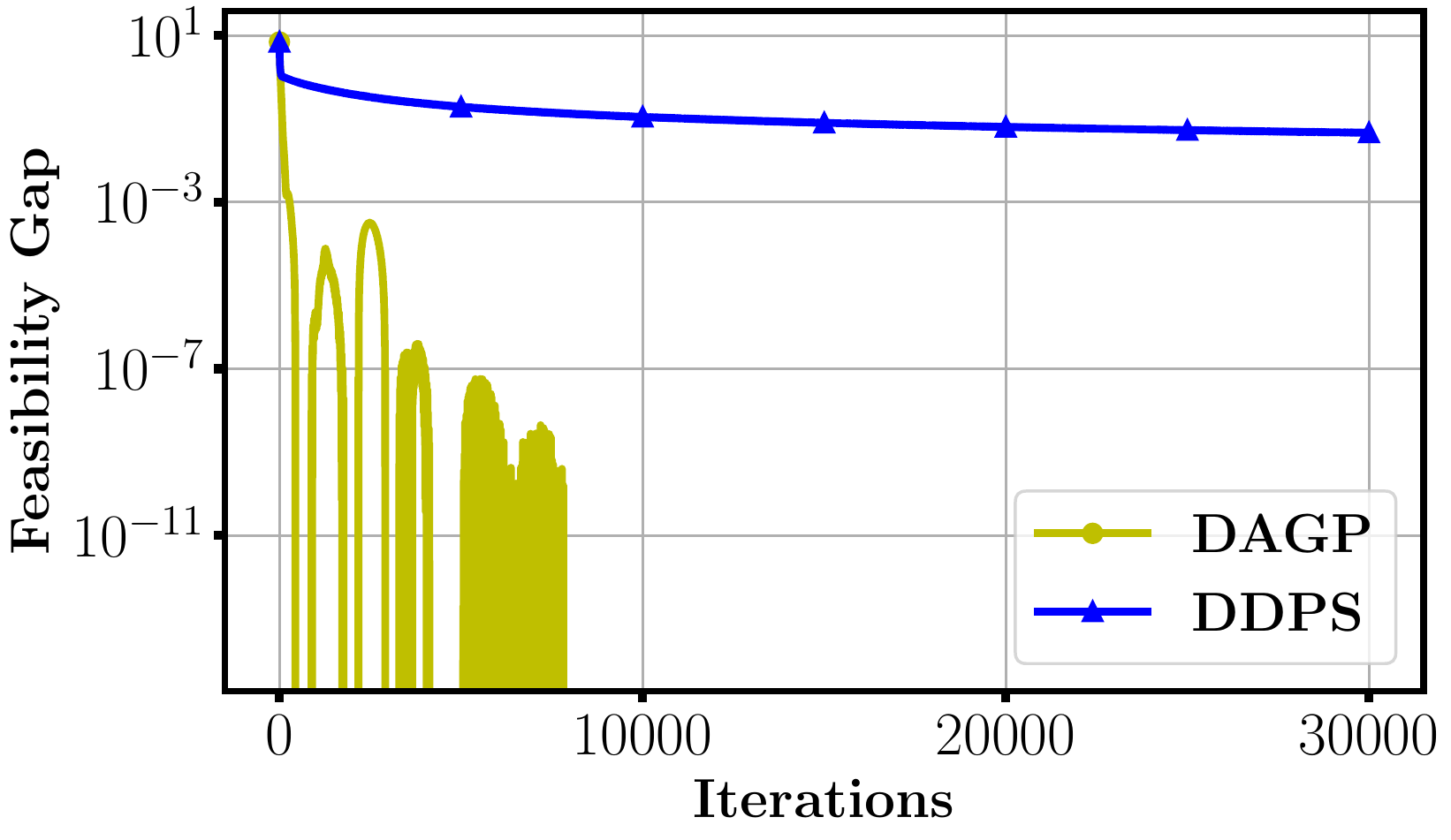} \label{fig: fesgpN20} }
    \caption{Second setup results with $m=10$ and $M=20$. As DDPS has not converged to a point in the feasible set, it can achieve less function value.}
    \label{fig:toy2}
\end{figure}

The results of these setups are shown respectively in Figures~\ref{fig:toy1} and \ref{fig:toy2}. We observe that
$\bar\bx$ moves completely to the feasible set in our algorithm unlike DDPS in which $\bar\bx$ becomes only close to the feasible set. 
Moreover, our algorithm converges faster to the optimal consensus solution in comparison with DDPS algorithm since DDPS needs a diminishing step size. 
To show that all nodes achieve consensus, the squared norm of error between $\bx^v$  at five random nodes and $\bx^0$ is plotted in Fig~\ref{fig: cons}, where all nodes converge to one stopping point. As described in section~\ref{sec: consensus_optimali}, $\sum_{ v \in \calV } \bg^v  $  should become equal to zero to have optimal solution. For this reason, the norm of this variable is computed and shown in Fig~\ref{fig: gsum}, which approaches zero as the algorithm proceeds. 

In the first  setup, $\bC = [\bc_i^T] \in \bbR^{M\times m}$  has a null space, and the feasible set of its corresponding optimization problem is larger than the feasible set of the second experiment. The linear convergence rate of DAGP in Fig~\ref{fig: objectiveN10}, can be due to a situation where the constraints are not active at the solution, and the algorithm attains the optimal value of the unconstrained version of the problem. Then, the overall unconstrained objective function may be strongly convex, explaining faster convergence. On the other hand, in the second experiment, some constraints are active, and the algorithm slows down in converging to a consensus and an optimal solution. In Fig~\ref{fig: ObjectiveN20}, DDPS  achieves a smaller objective value as its local variables remain infeasible. 

\subsection{Logistic Regression}

For the sake of comparison with other algorithms, we examine an unconstrained logistic regression problem. We consider the  MNIST 
 \cite{lecun-mnisthandwrittendigit-2010} 
dataset restricted to two digits to form a binary classification problem. We once again consider a random directed graph with $M=20$ nodes as shown in Fig~\ref{graph3}. 
Total number of $N_s=10000$ images are used for training the model, i.e. for minimizing the logistic loss function defined  as
\begin{equation}
\label{eq: lr}
    \minl_{\bw\in\bbR^{784}} \;\; \suml_{i = 1}^{N_s} \log \left( 1 + \exp{(-y_i\bx_i^T\bw)} \right) + \frac{\lambda}{2}\nt{\bw},
\end{equation}
where $\left\{ \bx_i, y_i \right\}_{i=1}^{N_s} \subseteq \bbR^{784} \times \{+1, -1\} $ is the set of training samples, and $\lambda$ is the regularization parameter, which is chosen to be $\slfrac{1}{N_s}$. We assume that the training samples are distributed in a balanced way between $20$ nodes; as a result, each node has $500$ training samples. The loss function at each node $f_v$ is the collection of terms in \eqref{eq: lr} associated with the samples of the node $v$.
The regularization term ensures that this problem is strongly convex, leading to linear rate of convergence for all algorithms.

\begin{figure}
    \centering
    \includegraphics[width = 0.75\linewidth]{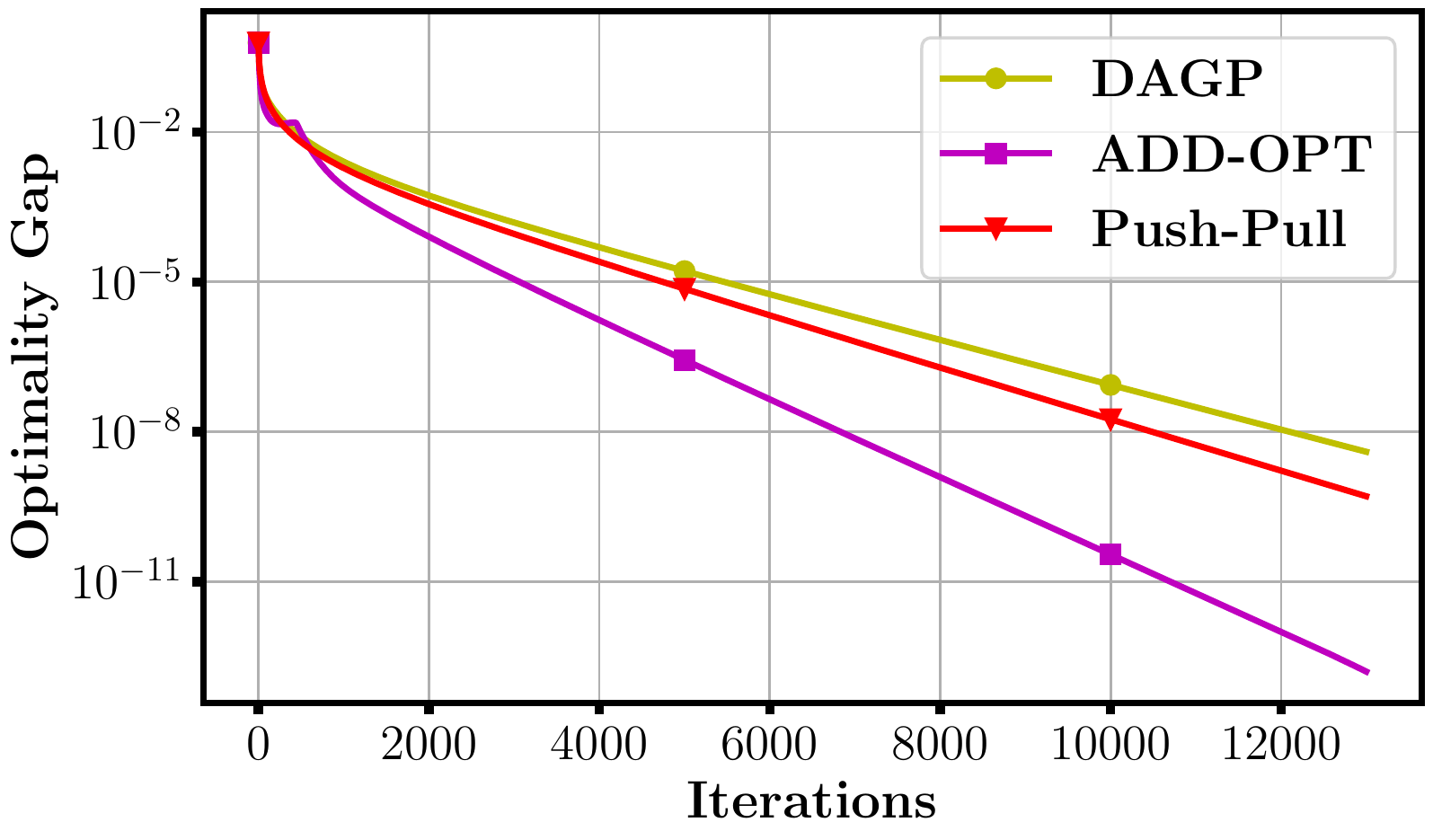}
    \caption{Convergence rate comparison of decentralized unconstrained algorithms over directed graphs. A fixed step size is used in all algorithms. }
    \label{fig: LR}
\end{figure}

We compare DAGP with Push-Pull and ADD-OPT. For all algorithms, fixed, similar, and  appropriate  step size is used. Centralized gradient descent is used to determine the optimal value $f^*$ and compare all algorithms objective values with respect to it. The results for optimality gap, defined as $\sum_{v\in\calV} f_v(\bar\bx) - f^*$, are shown in Fig~\ref{fig: LR}. 

As observed, the difference between the convergence rate of these three graphs is minimal.
As discussed earlier, optimal step sizes are not used in this experiment due to  the size of the problem. Nevertheless, for smaller experiments, we observe that DAGP and Push-Pull can utilize larger step sizes, while ADD-OPT fails to converge with a similar step size. We  conclude that by choosing the optimal values of the step size, DAGP and Push-Pull behave similarly, and both outperform ADD-OPT.
The practical applicability of our algorithm is immediate, as it is competitive with respect to constrained problems, but also provably capable of solving unconstrained optimization problems, without any modifications to the algorithmic structure.

\section{Conclusion}
\label{sec: conclusion}
We introduce the Double Averaging Gradient Projection (DAGP) algorithm designed for solving decentralized optimization problems over directed graphs for problems with and without constraints. 
In contrast to the existing literature, DAGP allows for different constraints at each node in a directed graph. Previous algorithms such as DDPS~\cite{xi2016distributed}, require a decreasing step size to solve similar constrained problems, while DAGP requires constant step size by employing a gradient tracking technique. 
By taking advantage of the projection operator  onto convex sets and double averaging, we prove the $\calO\left(\slfrac{1}{\sqrt{n}}\right)$ rate of convergence on constrained problems in a convex and smooth setting. Comparing to the previously proposed algorithms, DAGP is competitive with other decentralized unconstrained directed optimization algorithms, such as ADP-OPT~\cite{Xi2018} and Push-Pull~\cite{pu2020push}, and greatly outcompetes the previous decentralized constrained optimization algorithms such as DDPS.

\section{acknowledgement}
This work was partially supported by the Wallenberg AI, Autonomous Systems and Software Program (WASP) funded by the Knut and Alice Wallenberg Foundation.

\bibliographystyle{plain}
\bibliography{Arxiv}

\newpage
\appendix

The notations used throughout the paper are used similarly here, except for the $\bz^v$ variable, to which the iteration number is added.
\section*{Appendix A: Proof of Theorem 1}
We start by presenting the following lemma.
\begin{lemma} Let $\ker(\bW^T) = \ker(\bQ)$. Then, $\bQ\bW\bx = \bzero$ if and only if $\bx \in \ker(\bW)$.
\label{lemma: kernel} 
\begin{proof}
The forward proof is trivial. For the backward proof, we can write
\begin{align*}
    \bQ\bW\bx = \bzero & \Rightarrow \bW\bx \in \ker(\bQ) \\
     & \Rightarrow \bW\bx \in \ker(\bW^T) \\
     & \Rightarrow \bW^T\bW\bx = \bzero  \\
     & \Rightarrow \bx^T\bW^T\bW\bx = 0  \\
     & \Rightarrow \nt{\bW\bx}^2 = 0 \\
     & \Rightarrow \bx \in \ker(\bW)
\end{align*}
\end{proof}
\end{lemma}

Now, consider an arbitrary stopping point of the algorithm, that is $\bx_{n+1}^v = \bx^v_n = \bx^v $, $\nabla \mbf_{n+1} = \nabla \mbf_n = \nabla \mbf$, $ \bh^v_{n+1} = \bh^v_n = \bh^v $ and $ \bg^v_{n+1} = \bg^v_n = \bg^v$. 
We have
%
\begin{gather}
    \bZ     = \bX - \bW\bX - \mu \left( \nabla \mbf -  \bG   \right) \label{eq: z2} \\
    \bX     = \mathbfcal{P}_{S} \left( \bZ \right) \label{eq: x2} \\
    \rho\left[\nabla\mbf-\bG+\frac{1}{\mu}\left(\bZ-\bX\right)\right]+\alpha \left(\bH-\bG \right) = \textbf{O}  \label{eq: g2}\\
     \bQ\left(\bH - \bG\right) = \textbf{O}, \label{eq: h2}
\end{gather}
where $\bZ, \bG, \bH, \nabla\mbf, \mathbfcal{P}_{S}$ are matrices with $\bz^v, \bg^v, \bh^v,$ $ \nabla f^v(\bx^v), P_{S_v}(\bz^v)$ as their rows.
Left multiplying \eqref{eq: g2} by $\bQ$, considering \eqref{eq: h2}, we have
\begin{equation}
\label{eq: tmp1}
    \bQ(\bG-\nabla\mbf) = \frac{1}{\mu}\bQ(\bZ-\bX).
\end{equation}
Left multiplying \eqref{eq: z2} by $\bQ$, and applying \eqref{eq: tmp1} leads to $\bQ\bW\bX = \bzero$. Therefore, $\bX \in \ker(\bW)= \mathrm{span}\{\bone\}$ based on the result of Lemma \ref{lemma: kernel}, which means that $\bx^v = \bx^*,\; \forall v \in \calV$.
As $\bX \in \ker(\bW)$, \eqref{eq: z2} reduces to $\bZ - \bX  = \mu(\bG - \nabla\mbf)$, which leads to $\bH = \bG$ by incorporating it into \eqref{eq: g2}. Since \eqref{eq: h_update} is designed to preserve the summation of $\bh^v$s, and  each element of $\bH$ is initialized with zero vector, we have 
\begin{equation}
\label{eq: tmp2}
\bone^T\bG = \bone^T\bH = \suml_{v\in \calV} ({\bh^v})^T = \bzero^T. 
\end{equation}
From \eqref{eq: x2}, we have
$
\bZ-\bX \in \brond\bI_S $, consequently, $ \mu(\bG - \nabla\mbf) \in \brond\bI_S.$
%
As $\brond\bI_{S_v}$ is a cone, and therefore invariant to scaling, we can write $(\bG - \nabla\mbf) \in \brond\bI_S$. Left multiplying by $\bone^T$, and moving all the terms to one side, considering \eqref{eq: tmp2}, we have
\begin{equation}
    \bzero    \in \suml_{v \in \calV} \left( \partial I_{S_v}(\bx^*) + \nabla f_v(\bx^*) \right),
\end{equation} 
which shows that $\bx^*$ is the optimal solution. \qed

\newpage
\section*{Appendix B: Proof of Theorem 2}
We start by defining
\begin{equation}
    F^v(\bx)=f_v(\bx)-f_v(\bx^*)-\langle\nabla f_v(\bx^*),\bx-\bx^*\rangle
\end{equation}
and
\begin{equation}
    F^v_n=F^v(\bx^v_n).
\end{equation}
Note that from the convexity of $f_v$, the values of $F^v(\bx)$, particularly $F^v_n$ are non-negative. From convexity, we also conclude that
\begin{align}\label{eq:4}
    & F^v_n+\left\langle\nabla f_v(\bx^*)-\nabla f_v(\bx^v_n), \bx^v_n-\bx^*\right\rangle = \nonumber \\
    & f_v(\bx^v_n)-f_v(\bx^*)+\left\langle\nabla f_v(\bx^v_n), \bx^*-\bx^v_n\right\rangle\leq 0. 
\end{align}
%
From the $L-$smoothness property of $f_v$, we also obtain
\begin{align}
        & F_{n+1}^v-F_n^v-\langle\nabla f_v(\bx_n^v)-\nabla f_v(\bx^*),\bx_{n+1}^v-\bx^v_n\rangle = \nonumber \\
        & f_v(\bx^v_{n+1})-f_v(\bx^v_n)-\langle\nabla f_v(\bx_n^v),\bx_{n+1}^v-\bx^v_n\rangle \leq  \nonumber \\
        & \frac L2\left\|\bx_{n+1}^v-\bx_{n}^v\right\|^2. \label{eq:5}
\end{align}
Adding \eqref{eq:4} to \eqref{eq:5} yields:
\begin{equation}\label{eq:6}
    F^v_{n+1}+\left\langle\nabla f_v(\bx^*)-\nabla f_v(\bx^v_n),\ \bx^v_{n+1}-\bx^*\right\rangle 
    -\frac L2\left\|\bx_{n+1}^v-\bx_{n}^v\right\|^2\leq 0.
\end{equation}
Now, we define
\begin{equation}\label{eq:2}
    T^v(\bx)=-\langle\bn^v,\bx-\bx^*\rangle,
\end{equation}
and $T_n^v=T^v(\bx^v_{n}),$ where $\bn^v \in \partial I_{S_v}(\bx^*)$. The fact that $\bx_{n+1}^v\in S_v$ yields $T^v_{n+1}\geq 0$. Note that from $\bx^v_{n+1}=P_{S_v}(\bz^v_n)$ and the fact that $\bx^*\in S_v$, we have
\begin{equation}\label{eq:1}
    \langle\bx^*-\bx^v_{n+1},\bz_n^v-\bx^v_{n+1}\rangle\leq 0,
\end{equation}
which can also be written as
\begin{equation}\label{eq:3}
    \mu T^v_{n+1}+\langle\bx^*-\bx^v_{n+1},\bz_n^v-\bx^v_{n+1}-\mu\bn^v\rangle\leq 0.
\end{equation}

Multiplying \eqref{eq:6} by $\mu$, adding to \eqref{eq:3}, plugging the definition of $\bz_n^v$ in and summing over $v\in \calV$ and $n=0,1,\ldots,N-1$ we obtain
\begin{align}\label{eq:7}
    \mu & \suml_{n\in [N], v}\left(F_{n+1}^v+T_{n+1}^v\right)  -\frac {L\mu}2\suml_{n\in [N], v}\left\|\bx_{n+1}^v-\bx_{n}^v\right\|^2 \nonumber \\
    & + \suml_{n\in [N], v}\Big\langle\bx^*-\bx^v_{n+1},\bx_n^v-\suml_{u}w_{vu}\bx^u_n-\bx^v_{n+1}  +\mu(\bg^v_n-\nabla f_v(\bx^*)-\bn^v)\Big\rangle\leq 0.
\end{align}
We also replace the expression of $\bz^v_n$ in the dynamics of $\bg^v_n$ leading to
\begin{equation}\label{eq:8}
    \bg_{n+1}^v=\bg_n^v+\frac\rho\mu\left(\bx_n^v-\suml_{u}w_{vu}\bx^u_n-\bx^v_{n+1}\right)+\alpha\bdelta_n^v,
\end{equation}
where $\bdelta_n^v=\bh_n^v-\bg_n^v$ that follows the following dynamics
\begin{equation}\label{eq:9}
      \bdelta_{n+1}^v= (1-\alpha)\bdelta_{n}^v-\suml_uq_{vu}\bdelta^u_n 
     -\frac\rho\mu\left(\bx_n^v-\suml_{u}w_{vu}\bx^u_n-\bx^v_{n+1}\right).  
\end{equation}

For simplicity, we define $\tlbx_{n}^v=\bx^v_n-\bx^*$ and $\tlbg^v_n=\bg^v_n-\nabla f_v(\bx^*)-\bn^v$ and rewrite \eqref{eq:7}, \eqref{eq:8} and \eqref{eq:9} as

\begin{align}
    &\suml_{n=0}^{N-1}\left(\mu \suml_v\left(F_{n+1}^v+T_{n+1}^v\right)+\frac{\eta}2\suml_{u,v}\|\bx_{n+1}^u-\bx_{n+1}^v\|^2\right) \nonumber\\
    & - \suml_{n=0}^{N-1}\suml_v\left\langle\tlbx^v_{n+1},\tlbx_n^v-\suml_{u}w_{vu}\tlbx^u_n-\tlbx^v_{n+1}+\mu\tlbg^v_n\right\rangle \nonumber\\
    &  -\frac {L\mu}2 \suml_{n=0}^{N-1}\suml_v \left\|\tlbx_{n+1}^v-\tlbx_{n}^v\right\|^2 \nonumber \\
    & -\frac{\eta}2\suml_{n=0}^{N-1}\suml_{u,v}\|\tlbx_{n+1}^u-\tlbx_{n+1}^v\|^2 \leq 0, \label{eq:7p}
\end{align}

\begin{equation}\label{eq:8p}
    \tlbg_{n+1}^v=\tlbg_n^v+\frac\rho\mu\left(\tlbx_n^v-\suml_{u}w_{vu}\tlbx^u_n-\tlbx^v_{n+1}\right)+\alpha\bdelta_n^v,
\end{equation}

\begin{equation}
    \bdelta_{n+1}^v= (1-\alpha)\bdelta_{n}^v-\suml_uq_{vu}\bdelta^u_n  -\frac\rho\mu\left(\tlbx_n^v-\suml_{u}w_{vu}\tlbx^u_n-\tlbx^v_{n+1}\right), \label{eq:9p}
\end{equation}
 where in the first inequality we also add and remove the term $\frac{\eta}2\suml_{u,v}\|\bx_{n+1}^u-\bx_{n+1}^v\|^2 $. 
 
 In the following, we will consider the last three summations in \eqref{eq:7p} as  $A_N$. We show an asymptotic lower bound for $A_N$, i.e. show that there exists a constant $C$ only depending on the initial values such that for sufficiently large $N$, $A_N\geq -C$. Note that since $A_N\leq 0$, we must have $C\geq 0$. Then, we conclude from \eqref{eq:7p} that
\begin{equation}
    \suml_{n=0}^{N-1}\Bigg(\mu \suml_v\left(F_{n+1}^v+T_{n+1}^v\right)
     +\frac{\eta}2\suml_{u,v}\|\bx_{n+1}^u-\bx_{n+1}^v\|^2\Bigg)\leq C.
\end{equation}

Defining $\barbx^v_N=\frac 1N\suml_{n=0}^{N-1}\bx^v_n$ and noting that each term in the summation over $n$ is a fixed convex function of $\{\bx_{n+1}^v\}_v$, we may recall Jensen's inequality to conclude
 \begin{equation}
 \mu \suml_vF^v(\barbx^v_N)+T^v(\barbx^v_N)
 +\frac{\eta}2\suml_{u,v}\|\barbx_{N}^u-\barbx_{N}^v\|^2
    \leq \frac CN.
\end{equation}
Defining $\barbx_N=\frac 1M\suml_v\barbx^v_N$, we conclude that
\[
\|\barbx_N-\barbx^u_N\|^2=O(\frac 1N). 
\]
Since $\barbx^u_N\in S_u$, we also conclude that 
\[
\mathrm{dist}^2(\barbx_N, S_u)=O(\frac 1N). 
\]
Finally,
\begin{align}
    \left|\suml_vf_v(\barbx^v_N)-\suml_vf_v(\bx^*)\right| & \leq \frac{C}{\mu N}+\suml_v\left|\langle\bn^v+\nabla f_v(\bx^*),\barbx_N^v-\barbx_N\rangle\right| \nonumber \nwl 
    & \leq\frac{C}{\mu N}+\sqrt{\suml_v\|\bn^v+\nabla f_v(\bx^*)\|^2}\sqrt{\suml_v\|\barbx_N^v-\barbx_N\|^2} \nonumber \nwl 
    & =O(\frac 1{\sqrt{N}}) \nonumber
\end{align}
which completes the proof. \qed

\section*{Bound on $A_N$}
To find the bound $C$, we start by simplifying the notation in \eqref{eq:7p}, \eqref{eq:8p} and \eqref{eq:9p}. Let us introduce 
\begin{equation}
    \bPsi_n= \left[
\begin{array}{cccc} 
\tlbX_{n+1} & \tlbX_{n}&\tlbG_{n} & \bDelta_{n}
\end{array}
\right]^T,
\end{equation}
where $\tlbX_{n}, \tlbG_{n}, \bDelta_{n}$ are matrices with $\tlbx^v_{n}, \tlbg^v_{n}, \bdelta^v_{n}$ as their $v^\tth$ row, respectively. We may write \eqref{eq:8p} and \eqref{eq:9p} as
\begin{equation}\label{eq:mat}
    \bPsi_{n+1}=\bR\bPsi_n+\bP\tlbX_{n+2},  \qquad  n = 0, \dots, N-2 
\end{equation} 
where $\bP$ and $\bR$ are defined as
\begin{equation} \label{eq: R_define}
\bR = \left[
\begin{array}{cccc}
        \bO & \bO & \bO & \bO   \\
        \bI & \bO & \bO & \bO  \\
        -\frac{\rho}{\mu}\bI & \frac{\rho}{\mu}(\bI-\bW) & \bI & \alpha \bI \\
        \frac{\rho}{\mu}\bI & -\frac{\rho}{\mu}(\bI-\bW) & \bO & (1-\alpha)\bI-\bQ 
\end{array} \right],
\qquad
\bP = \left[
\begin{array}{c}
     \bI  \\
     \bO  \\
     \bO  \\
     \bO 
\end{array}
\right].
\end{equation} 
We also have
\begin{equation}
    A_N=\suml_{n=0}^{N-1}\left\langle\bPsi_n,\bS\bPsi_n\right\rangle,
\end{equation}
where $\bS$ is computed as
\begin{equation}\label{eq: S_define}
    \bS=\left[\begin{array}{cccc}
        \left(1-\frac{L\mu}2\right)\bI-M\eta\left(\bI-\frac 1M\bone\bone^T\right) & -\frac 12(\bI-\bW)+\frac {L\mu}2\bI  & -\frac {\mu}2\bI  & \bO\\
        -\frac 12(\bI-\bW^T)+\frac {L\mu}2\bI & -\frac{L\mu}2\bI & \bO &\bO \\
        -\frac{\mu}2\bI & \bO & \bO & \bO \\
         \bO & \bO & \bO& \bO
    \end{array}\right].
\end{equation}

The following Lemmas provide a  guarantee that $A_N$ is bounded.

\begin{lemma}\label{lem1}
Consider matrices $\bR, \bP$ and $\bS$,  defined in 
(\ref{eq: R_define}, \ref{eq: S_define}). Define a "dual" sequence $\{\bLambda_n\}_{n=-1}^{N-1}$ such that $\bLambda_{N-1}= \bLambda_{-1} = \bO$.
Suppose that there exists a $C\geq 0$ such that
for every $\beta>0$, the system of equations in  \eqref{eq:mat} together with
%
%
\begin{align}
    \bLambda_{n-1}-\bR^T\bLambda_n+(\bS+\left(C+\beta\right)\delta_{n,0}\bI)\bPsi_{n}=\bO & \qquad n=0,1,\ldots,N-1 \label{eq:10}\\
    \bP^T\bLambda_n + \beta\tlbX_{n+2} = \bO \qquad\qquad\quad & \qquad n=0,1,\ldots, N-2 \label{eq:13}
\end{align}

has no non-zero solution for $\{\bPsi_n,\bLambda_n,\tlbX_{n+2}\}$.
Then, $\bA_n\geq -C\|\bPsi_0\|_{\mathrm{F}}^2$ always holds true.
\begin{proof}
Note that the claim is equivalent to the statement that zero is the optimal value for the optimization problem
\begin{align}
    \minl_{\{\bPsi_n\}_{n=0}^{N-1},\{\tlbX_{n+2}\}_{n=0}^{N-2}} & \quad \frac{1}{2}\suml_{n=0}^{N-1}\langle\bPsi_n,\bS\bPsi_n\rangle+\frac{C}{2}\|\bPsi_0\|^2_{\mathrm{F}} \nonumber\\
   \text{subject to}\qquad\; & \quad \bPsi_{n+1}=\bR\bPsi_n+\bP\tlbX_{n+2},\quad n=0,1,\ldots, N-2
\end{align}
If the claim does not hold, the optimization is unbounded and the following restricted optimization will achieve a strictly negative optimal value at a non-zero solution:
\begin{align}
    \minl_{\{\bPsi_n\}_{n=0}^{N-1},\{\tlbX_{n+2}\}_{n=0}^{N-2}} & \quad \frac{1}{2}\suml_{n=0}^{N-1}\langle\bPsi_n,\bS\bPsi_n\rangle+\frac{C}{2}\|\bPsi_0\|^2_{\mathrm{F}} \nonumber \\
    \text{subject to} \qquad\; & \quad \bPsi_{n+1}=\bR\bPsi_n+\bP\tlbX_{n+2},\quad n=0,1,\ldots, N-2 \nonumber \\
    & \quad \frac{1}{2}\nf{\bPsi_0}^2 + \frac{1}{2}\suml_{n=0}^{N-2}\|\tlbX_{n+2}\|_\mathrm{F}^2 \leq \frac{1}{2}
\end{align}
Such a solution satisfies the KKT condition,  which coincides with \eqref{eq:10}, \eqref{eq:13} where $\{\bLambda_n\},\beta\geq 0$ are dual (Lagrangian) multipliers corresponding to the constraints. We also observe that the optimal value at this point is given by $-\beta\left( \nf{\bPsi_0}^2 + \suml_{n=0}^{N-2}\|\tlbX_{n+2}\|_\mathrm{F}^2\right)$. This shows that $\beta>0$. This contradicts the assumption that such a point does not exists and completes the proof. 
\end{proof}
\end{lemma}
We may further simplify the conditions in Lemma~\ref{lem1} by the following result:
\begin{lemma}
For a complex value $z$ and a real value $\beta$, define
\begin{equation}
    \bF(z,\beta)=\left[\begin{array}{ccc}
        \bS & z^{-1}\bI-\bR^T & \bO\\
        z\bI-\bR & \bO & -\bP \\
        \bO & -\bP^T & -\beta\bI
    \end{array}\right].
\end{equation}
The condition of Lemma~\ref{lem1} is satisfied if the matrix
\begin{equation}
    \liml_{z\to 0}\left[\begin{array}{ccc}
        \bI\ & \bO\ & \bO   
    \end{array}\right]\bF^{-1}\left[\begin{array}{c}
          -(C+\beta)\bI \\
            \bI\\
            \bO
    \end{array}\right]
\end{equation}
doesn't have an eigenvalue of 1.
\begin{proof}
With an abuse of notation, define the z-transforms
\begin{equation*}
    \bPsi(z)=\suml_{n=0}^{N-1}\bPsi_nz^{-n},\quad \bLambda(z)=\suml_{n=0}^{N-2}\bLambda_nz^{-n}, \quad \bU(z)=\suml_{n=0}^{N-2}\tlbX_{n+2}z^{-n}
\end{equation*}
Then, for the sequences defined in (\ref{eq:10},~\ref{eq:13},~\ref{eq:mat}), we have
\begin{equation}
    (z^{-1}\bI-\bR^T)\bLambda(z)+\bS\bPsi(z)+(C+\beta)\bPsi_0=\bO
\end{equation}
\begin{equation}
    (z\bI-\bR)\bPsi(z)-\bPsi_0+\bR\bPsi_{N-1}z^{n-1}-\bP\bU(z)=\bO
\end{equation}
\begin{equation}
    \bP^T\bLambda(z) + \beta\bU(z) = \bO
\end{equation}
which can also be written as 
\begin{equation}
    \bF(z,\beta)\left[\begin{array}{c}
        \bPsi(z)\\
        \bLambda(z)\\
        \bU(z)
    \end{array}\right]=\left[\begin{array}{c}
        -(C+\beta)\bPsi_0\\
        \bPsi_0-\bR\bPsi_{N-1}z^{N-1}\\
        \bO
    \end{array}\right]
\end{equation}
Note that $\bF(z,\beta)$ may be rank-deficient at a finite number of points. Hence, there exists a sufficiently small simple loop $\calC$ around $z=0$ such that $\bF$ is invertible on and inside it except at $z=0$. In this region we have
\begin{equation}\label{eq:f}
    \bPsi(z)=\left[\begin{array}{ccc}
        \bI\ & \bO\ & \bO   
    \end{array}\right]\bF^{-1}\bA,
\end{equation}
where $\bA$ is 
\begin{equation*}
    \left[\begin{array}{c}
          -(C+\beta)\bI \\
            \bI\\
            \bO
    \end{array}\right]\bPsi_0+\left[\begin{array}{c}
          \bO \\
            -\bR\\
            \bO
    \end{array}\right]z^{N-1}\bPsi_{N-1}
\end{equation*}
On the other hand, 
\begin{equation}\label{eq:ff}
    2\pi j\bPsi_0=\ointl_{\calC}\frac 1z\bPsi(z)\td z
\end{equation}
Further from the Cauchy integral formula for sufficiently large $N$, we have
\begin{equation}\label{eq:fff}
    \ointl_{\calC}\frac 1z z^{N-1}\bF^{-1}(z,\beta)\td z=2\pi j\liml_{z\to 0}z^{N-1}\bF^{-1}(z,\beta)=\bO.
\end{equation}
By applying \eqref{eq:f} to \eqref{eq:ff}, considering relation in \eqref{eq:fff}, we can conclude that
\begin{equation}
    \bPsi_0=\liml_{z\to 0}\left[\begin{array}{ccc}
        \bI\ & \bO\ & \bO   
    \end{array}\right]\bF^{-1}\left[\begin{array}{c}
          -(C+\beta)\bI \\
            \bI\\
            \bO
    \end{array}\right]\bPsi_0
\end{equation}
Note that by the assumption we can conclude that $\bPsi_0=\bO$, which completes the proof.
\end{proof}
\end{lemma}

\end{document}